\documentclass[12pt]{article}

\usepackage{amsmath}
\usepackage{amsfonts}

\newtheorem{Theorem}{Theorem}[section]

\newtheorem{Lemma}{Lemma}[section]

\newenvironment {Proof} {\noindent {\bf Proof.}}{\hspace*{\fill}$\Box$\par\vspace{4mm}}

\allowdisplaybreaks[4]
\begin{document}

\title{\bf  On Sum--Connectivity Index of Bicyclic Graphs}
\author{Zhibin Du, Bo Zhou\footnote{Corresponding author.
E-mail address: zhoubo@scnu.edu.cn}\\
\vspace{1mm}\\
Department of Mathematics, South China Normal University, \\
Guangzhou 510631, P. R. China }

\date{}
\maketitle

\begin{abstract}
We determine the minimum sum--connectivity index of bicyclic graphs
with $n$ vertices and matching number $m$, where $2\le m\le
\lfloor\frac{n}{2}\rfloor$, the minimum and the second minimum, as
well as  the maximum and the second maximum sum--connectivity
indices of bicyclic graphs with $n\ge 5$ vertices. The  extremal
graphs are characterized.\\ \\
{\em MSC 2000:} 05C90; 05C35; 05C07\\ \\
{\em Keywords:\/}  Sum--connectivity index; Randi\'{c} connectivity index; Matching number;
Bicyclic graphs
\end{abstract}

\baselineskip=0.24in

\section{Introduction}

The Randi\'{c} connectivity index \cite{Ran} is one of the most
successful molecular descriptors in structure--property and
structure--activity relationships studies, e.g., \cite{Ran2,ToCo}.
Its mathematical properties as well as those of its generalizations
have  been studied extensively as summarized in the books
\cite{LiGu,GuFu}.  Recently, a closely related variant of Randi\'{c}
connectivity index called the sum--connectivity index was proposed
in \cite{ZhTr}.

Let $G$ be a simple graph with vertex set $V(G)$   and edge set
$E(G)$. For $u\in V(G)$, $d_G(u)$  denotes the degree of $u$ in $G$.
The Randi\'{c} connectivity index (or product--connectivity index \cite{ZhTr,LTZ})
of the graph $G$ is defined as
\cite{Ran}
\[
R(G)=\sum_{uv\in E(G)}\frac{1}{\sqrt{d_G(u)d_G(v)}}\,.
\]
The sum--connectivity index of $G$ is defined as \cite{ZhTr}
\[
\chi (G)=\sum_{uv\in E(G)}\frac{1}{\sqrt{d_G(u)+d_G(v)}}\,.
\]

It has been found that the sum--connectivity index and the
Randi\'{c} connectivity index correlate well among themselves and
with $\pi$-electronic energy of benzenoid hydrocarbons \cite{ZhTr,LTZ}.
Some mathematical properties of the sum--connectivity index have
been established in \cite{ZhTr,DZ1,DZ2}. Recall that an $n$-vertex
connected graph is known as a tree, a unicyclic graph and a bicyclic
graph if it possesses $n-1$, $n$ and $n+1$ edges, respectively. In \cite{DZ1}, we
obtained the minimum sum--connectivity indices of trees and
unicyclic graphs respectively with given number of vertices and
matching number,
 and determined the corresponding extremal graphs. The
$n$-vertex trees with the first a few minimum and maximum
sum--connectivity indices were determined in \cite{ZhTr}, while the
$n$-vertex unicyclic graphs with the first a few minimum and maximum
sum--connectivity indices were determined in \cite{DZ1} and
\cite{DZ2}, respectively. In this paper, we consider the
sum--connectivity indices of bicyclic graphs.

A matching $M$ of the graph $G$ is a subset of $E(G)$ such that no
two edges in $M$ share a common vertex. A matching $M$ of $G$ is
said to be maximum, if for any other matching $M'$ of $G$, $|M'|\le
|M|$. The matching number of $G$ is the number of edges of a maximum
matching in $G$.

If $M$ is a matching of a graph $G$ and vertex $v\in V(G)$ is
incident with an edge of $M$, then $v$ is said to be $M$-saturated,
and if every vertex of $G$ is $M$-saturated, then $M$ is a perfect
matching.

In this paper, we obtain the minimum sum--connectivity index in  the
set of bicyclic graphs with $n$ vertices and matching number $m$,
where $2\le m \le \lfloor n/2\rfloor$. We also determine the minimum
and the second minimum, as well as the maximum and the second
maximum sum--connectivity indices in the set of bicyclic graphs with
$n\ge 5$ vertices. The extremal graphs are characterized.

Study on the Randi\'{c} connectivity indices of bicyclic graphs may be
found in \cite{LiGu,ZLW,CGHP,WZL}, and in particular, the minimum
and the maximum Randi\'{c} connectivity indices in the set of
bicyclic graphs with $n\ge 5$ vertices were determined in \cite{WZL} and \cite{CGHP}, respectively.

We note that some other graph invariants based on end--vertex
degrees of edges in a graph have been studied recently,
see, e.g., \cite{FGV,VF, YZT}.

\section{Preliminaries}

For $2\le m \le \lfloor n/2\rfloor$, let $\mathcal {B}(n,m)$ be the
set of bicyclic graphs with $n$ vertices and matching number $m$.

For $3\le m \le \lfloor n/2\rfloor$, let $B_{n,m}$ be the graph
obtained by identifying a vertex of two triangles, and attaching
$n-2m+1$ pendent vertices (vertices of degree one) and $m-3$ paths
on two vertices to the common vertex of the two triangles, see Fig.
1. Obviously, $B_{n,m}\in\mathcal {B}(n,m)$.

\begin{center}
\setlength{\unitlength}{0.75mm} \thicklines
{\begin{picture}(29,40)(-8,0) \put(0,0){\begin{picture}(110,50)


\put(10,22){\circle*{2.5}}
\put(-6.5,31.5){\circle*{2.5}}\put(26.5,31.5){\circle*{2.5}}
\put(-6.5,12){\circle*{2.5}}\put(26.5,12){\circle*{2.5}}
\put(10,22){\line(-5,3){16}}\put(10,22){\line(-5,-3){16}}
\put(10,22){\line(5,3){16}}\put(10,22){\line(5,-3){16}}
\put(10,22){\line(3,-4){8}}\put(10,22){\line(-3,-4){8}}
\put(-6.5,31.5){\line(-1,0){15}}\put(-6.5,12){\line(-1,0){15}}
\put(-22,12){\circle*{2.5}}\put(-22,31.5){\circle*{2.5}}
\put(2.5,12){\circle*{2.5}}\put(17.5,12){\circle*{2.5}}\put(17,12){\line(-1,0){15}}
\put(-22,18){\circle*{0.5}}\put(-22,21.5){\circle*{0.5}}\put(-22,25){\circle*{0.5}}
\put(26.5,18){\circle*{0.5}}\put(26.5,21.5){\circle*{0.5}}\put(26.5,25){\circle*{0.5}}
\put(29,20){$\left\} \begin{array}{c} \\ \\ \\
\end{array} \right.$} \put(35,20){$n-2m+1$}
\put(-29,20){$\left\{ \begin{array}{c} \\ \\ \\
\end{array} \right.$} \put(-44,20){$m-3$}
\put(-22,-10){Fig. 1. The graph $B_{n,m}$.}
\put(10,22){\line(-3,4){7.5}}\put(10,22){\line(3,4){7.5}}
\put(3,31.5){\circle*{2.5}}\put(17,31.5){\circle*{2.5}}
\put(3,31.5){\line(6,0){15}}
\end{picture}}
\end{picture}}
\end{center}

\vspace{1cm}
Let $C_n$ ba a cycle on $n\ge 3$ vertices. 
Let $\widetilde{\mathbb{B}}(n)$ be the set of $n$-vertex bicyclic
graphs without pendent vertices, where $n\ge 4$. Let ${\bf
B}_1^{(1)}(n)$ be the set of bicyclic graphs obtained by joining two
vertex--disjoint cycles $C_a$ and $C_b$ with $a+b=n$ by an edge,
where $n\ge 6$. Let ${\bf B}_1^{(2)}(n)$ be the set of bicyclic
graphs obtained by joining two vertex--disjoint cycles $C_a$ and
$C_b$ with $a+b<n$ by a path of length $n-a-b+1$, where $n\ge7$. Let
${\bf B}_2(n)$ be the set of bicyclic graphs obtained by identifying
a vertex of $C_a$ and a vertex of $C_b$ with $a+b=n+1$, where
$n\ge5$. Let ${\bf B}_3^{(1)}(n)$ be the set of bicyclic graphs
obtained from $C_n$ by adding an edge, where $n\ge 4$. Let ${\bf
B}_3^{(2)}(n)$ be the set of bicyclic graphs obtained by joining two
non-adjacent vertices of $C_a$ with $4\le a \le n-1$ by a path of
length $n-a+1$, where $n\ge 5$. Obviously,
$\widetilde{\mathbb{B}}(n) ={\bf B}_1^{(1)}(n)\cup {\bf
B}_1^{(2)}(n)\cup {\bf B}_2(n)\cup {\bf B}_3^{(1)}(n)\cup {\bf
B}_3^{(2)}(n)$.

Let $\mathbb{B}(n)$ be the set of bicyclic graphs on $n\ge 4$
vertices.

\section{Minimum sum--connectivity index of bicyclic graphs with given matching number}

First we give some lemmas that will be used.

For a graph $G$ with $u\in V(G)$, $G-u$ denotes the graph resulting
from $G$ by deleting the vertex $u$ (and its incident edges).

\begin{Lemma} \cite{DZ1} \label{l5}
Let $G$ be an $n$-vertex connected graph with a pendent vertex $u$,
where $n\ge 4$. Let $v$ be the unique neighbor of $u$, and let $w$
be a neighbor of $v$ different from $u$.
\begin{enumerate}

\item[(i)]
If $d_G(v)=2$ and there is at most one pendent neighbor of $w$ in
$G$, then
\[
\chi(G)-\chi(G-u-v)\ge
\frac{d_G(w)-1}{\sqrt{d_G(w)+2}}-\frac{d_G(w)-3}{\sqrt{d_G(w)+1}}-\frac{1}{\sqrt{d_G(w)}}+\frac{1}{\sqrt{3}}
\]
with equality if and only if one neighbor of $w$ has degree one, and
the other neighbors of $w$ are of degree two.

\item[(ii)]
If there are at most $k$ pendent neighbors of $v$ in $G$, then
\[
\chi(G)-\chi(G-u)\ge
\frac{d_G(v)-k}{\sqrt{d_G(v)+2}}+\frac{2k-d_G(v)}{\sqrt{d_G(v)+1}}-\frac{k-1}{\sqrt{d_G(v)}}
\]
with equality if and only if $k$ neighbors of $v$ have degree one,
and the other neighbors of $v$ are of degree two.
\end{enumerate}
\end{Lemma}

\begin{Lemma} \cite{DZ1} \label{l6}
(i) The function
$\frac{x-1}{\sqrt{x+2}}-\frac{x-3}{\sqrt{x+1}}-\frac{1}{\sqrt{x}}$
is decreasing for $x\ge2$.

(ii) For integer $a\ge 1$, the function
$\frac{x-a}{\sqrt{x+2}}+\frac{2a-x}{\sqrt{x+1}}-\frac{a-1}{\sqrt{x}}$
is decreasing for $x\ge a+1$.
\end{Lemma}

\begin{Lemma} \cite{DZ1} \label{l7}
Let $G$ be a connected graph with $uv\in E(G)$, where $d_G(u)$,
$d_G(v)\ge2$, and $u$ and $v$ have no common neighbor in $G$. Let
$G_1$ be the graph obtained from $G$ by deleting the edge $uv$,
identifying $u$ and $v$, which is denoted by $w$, and attaching a
pendent vertex to $w$. Then $\chi(G)>\chi(G_1)$.
\end{Lemma}

\begin{Lemma}\label{l11}
For $m\ge3$,
$m+\frac{4}{\sqrt{6}}-\frac{3}{2}>\frac{m+1}{\sqrt{m+4}}+\frac{1}{\sqrt{m+3}}+\frac{m-3}{\sqrt{3}}+1$,
and for $m\ge5$,
$(\frac{1}{2}+\frac{1}{\sqrt{6}})m-\frac{1}{2}-\frac{2}{\sqrt{6}}+\sqrt{2}
>\frac{m+1}{\sqrt{m+4}}+\frac{1}{\sqrt{m+3}}+\frac{m-3}{\sqrt{3}}+1$.
\end{Lemma}

\begin{Proof}
Let $f(m)=\left(m+\frac{4}{\sqrt{6}}-\frac{3}{2}\right)
-\left(\frac{m+1}{\sqrt{m+4}}+\frac{1}{\sqrt{m+3}}+\frac{m-3}{\sqrt{3}}+1\right)$
for $m\ge 3$, and let
$g(m)=\left[\left(\frac{1}{2}+\frac{1}{\sqrt{6}}\right)m-\frac{1}{2}-\frac{2}{\sqrt{6}}+\sqrt{2}\right]
-\left(\frac{m+1}{\sqrt{m+4}}+\frac{1}{\sqrt{m+3}}+\frac{m-3}{\sqrt{3}}+1\right)$
for $m\ge5$. Note that
$f''(m)=g''(m)=-\frac{3}{4}(m+3)^{-5/2}+(\frac{1}{4}m+\frac{13}{4})(m+4)^{-5/2}>0$.
Then $f'(m)\ge f'(3)>0$, implying that $f(m)\ge f(3)>0$, and
$g'(m)\ge g'(5)>0$, implying that $g(m)\ge g(5)>0$.
\end{Proof}

\begin{Lemma}\label{l12}
For $m\ge3$,
\[
-\frac{m+1}{\sqrt{m+4}}+\frac{m-1}{\sqrt{m+3}}+\frac{1}{\sqrt{m+2}}
\ge-\frac{4}{\sqrt{7}}+\frac{2}{\sqrt{6}}+\frac{1}{\sqrt{5}}
\]
with equality if and only if $m=3$.
\end{Lemma}

\begin{Proof}
Let $f(m)=(m+2)^{-1/2}+m\cdot (m+3)^{-1/2}$ for $m\ge3$. Then
$f''(m)=\frac{3}{4}(m+2)^{-5/2}-(\frac{1}{4}m+3)(m+3)^{-5/2}<0$,
implying that $f(m)-f(m+1)$ is increasing on $m$. It is easily seen
that
\begin{eqnarray*}
&&-\frac{m+1}{\sqrt{m+4}}+\frac{m-1}{\sqrt{m+3}}+\frac{1}{\sqrt{m+2}}\\
&=&f(m)-f(m+1)\\
&\ge&f(3)-f(4)\\
&=&-\frac{4}{\sqrt{7}}+\frac{2}{\sqrt{6}}+\frac{1}{\sqrt{5}}
\end{eqnarray*}
with equality if and only if $m=3$.
\end{Proof}

Let $H_6$ be the graph obtained by attaching a pendent vertex to
every vertex of a triangle. For $2\le m \le \lfloor n/2\rfloor$, let
$U_{n,m}$ be the unicyclic graph obtained by attaching $n-2m+1$
pendent vertices and $m-2$ paths on two vertices to one vertex of a
triangle.

\begin{Lemma} \cite{DZ1} \label{l13} Let $G$ be a unicyclic graph
with $2m$ vertices and perfect matching, where $m\ge 3$. Suppose
that $G\ne H_6$. Then
\[
\chi(G)\ge\frac{m}{\sqrt{m+3}}+\frac{1}{\sqrt{m+2}}+\frac{m-2}{\sqrt{3}}+\frac{1}{2}
\]
with equality if and only if $G=U_{2m,m}$.
\end{Lemma}

For an edge $uv$ of the graph $G$ (the complement of $G$,
respectively), $G-uv$ ($G+uv$, respectively) denotes the graph
resulting from $G$ by deleting (adding, respectively) the edge $uv$.

\begin{Lemma}\label{l8}
Let $G\in\mathcal{B}(2m,m)$ and no pendent vertex has neighbor of
degree two, where  $m\ge 3$. Then
$\chi(G)\ge\frac{m+1}{\sqrt{m+4}}+\frac{1}{\sqrt{m+3}}+\frac{m-3}{\sqrt{3}}+1$
with equality if and only if $m=3$ and $G=B_{6,3}$.
\end{Lemma}

\begin{Proof} Let
$f(m)=\frac{m+1}{\sqrt{m+4}}+\frac{1}{\sqrt{m+3}}+\frac{m-3}{\sqrt{3}}+1$.

Since  $G\in\mathcal{B}(2m,m)$ and no pendent vertex has neighbor of
degree two,
$G$ is obtainable by attaching some pendent vertices
to a graph in $\widetilde{\mathbb{B}}(k)$, where $m\le k \le 2m$,
and any two pendent vertices have no common neighbor (if $k=2m$,
then no pendent vertex is attached).

\noindent {\bf Case 1.} There is no vertex of degree two in $G$.
Then either $k=m$, $G$ is obtainable by attaching a pendent
vertex to every vertex of a graph in $\widetilde{\mathbb{B}}(m)$, or
$k=m+1$,  $G$ is obtainable by attaching a pendent vertex to
every vertex with degree two of a graph in ${\bf B}_1^{(1)}(m+1)\cup
{\bf B}_3^{(1)}(m+1)$. By direct calculation, we find that  $\chi(G)=
\frac{5}{\sqrt{6}}+1>f(3)$ for $m=3$,  $\chi(G)\ge
\frac{1}{\sqrt{8}}+\frac{4}{\sqrt{7}}+\frac{2}{\sqrt{5}}+1> f(4)$
for $m=4$, and  $\chi(G)\ge
\left(\frac{1}{2}+\frac{1}{\sqrt{6}}\right)m-\frac{1}{2}-\frac{2}{\sqrt{6}}+\sqrt{2}$
for $m\ge 5$. Thus by Lemma \ref{l11}, we have $\chi(G)> f(m)$.

\noindent {\bf Case 2.} There is a vertex, say $u$, of degree two in
$G$. Denote by $v$ and $w$ the two neighbors of $u$ in $G$.  Then
one of the two edges incident with $u$, say $uv\in M$,  where $M$ is
a perfect matching of $G$.
Suppose that there is no vertex of degree two in any cycle of $G$.
Since no pendent vertex has neighbor of degree two in $G$,
$u$ lies on the path joining the two disjoint cycles of $G$.
For $G_1=G-uw+vw\in \mathcal{B}(2m,m)$,
the difference of the numbers of vertices of degree two outside
any cycle of $G$ and $G_1$ is equal to one, and thus by Lemma
\ref{l7}, $\chi(G_1)<\chi(G)$. Repeating the operation from $G$ to $G_1$,
we finally get a graph $G'\in \mathcal{B}(2m,m)$,
which has no vertex of degree two, such that $\chi(G)>\chi(G')$,
and thus the result follows from  Case 1.
Now suppose that $u$ lies on some cycle of $G$.
Consider $G'=G-uw$, which is a unicyclic graph with perfect
matching. If $G'=H_6$, then $G$ is obtained from $H_6$ by adding an
edge either between two pendent vertices, and thus
$\chi(G)=\frac{3}{\sqrt{6}}+\frac{2}{\sqrt{5}}+1$, or between two
neighbors of a vertex of degree three, one of which being a pendent
vertex, and thus
$\chi(G)=\frac{2}{\sqrt{7}}+\frac{2}{\sqrt{6}}+\frac{2}{\sqrt{5}}+\frac{1}{2}$.
In either case, $\chi(G)>f(3)$. Suppose that $G'\ne H_6$.  Then by
Lemma \ref{l13},
$\chi(G')\ge\frac{m}{\sqrt{m+3}}+\frac{1}{\sqrt{m+2}}+\frac{m-2}{\sqrt{3}}+\frac{1}{2}$.
Note that $2\le d_G(v),d_G(w)\le 5$ and $w$ has at most one pendent
neighbor. By Lemmas \ref{l6} (i) and \ref{l12}, we have
\begin{eqnarray*}
\chi(G)&=&\chi(G')+\frac{1}{\sqrt{d_G(w)+2}}+\left(\frac{1}{\sqrt{d_G(v)+2}}-\frac{1}{\sqrt{d_G(v)+1}}\right)\\
&&+\sum_{xw\in
E(G')}\left(\frac{1}{\sqrt{d_G(w)+d_G(x)}}-\frac{1}{\sqrt{d_G(w)+d_G(x)-1}}\right)\\
&\ge&\chi(G')+\frac{1}{\sqrt{d_G(w)+2}}+\left(\frac{1}{\sqrt{2+2}}-\frac{1}{\sqrt{2+1}}\right)\\
&&+\left[\frac{1}{\sqrt{d_G(w)+1}}-\frac{1}{\sqrt{d_G(w)+1-1}}\right.\\
&&\left.+(d_G(w)-2)\left(\frac{1}{\sqrt{d_G(w)+2}}-\frac{1}{\sqrt{d_G(w)+2-1}}\right)\right]\\
&=&\chi(G')+\left(\frac{d_G(w)-1}{\sqrt{d_G(w)+2}}-\frac{d_G(w)-3}{\sqrt{d_G(w)+1}}
-\frac{1}{\sqrt{d_G(w)}}\right)+\frac{1}{2}-\frac{1}{\sqrt{3}}\\
&\ge&\left(\frac{m}{\sqrt{m+3}}+\frac{1}{\sqrt{m+2}}+\frac{m-2}{\sqrt{3}}+\frac{1}{2}\right)\\
&&+\left(\frac{5-1}{\sqrt{5+2}}-\frac{5-3}{\sqrt{5+1}}-\frac{1}{\sqrt{5}}\right)+\frac{1}{2}-\frac{1}{\sqrt{3}}\\
&=&\frac{m}{\sqrt{m+3}}+\frac{1}{\sqrt{m+2}}+\frac{m-2}{\sqrt{3}}+1
-\frac{1}{\sqrt{3}}+\left(\frac{4}{\sqrt{7}}-\frac{2}{\sqrt{6}}-\frac{1}{\sqrt{5}}\right)\\
&\ge&\frac{m}{\sqrt{m+3}}+\frac{1}{\sqrt{m+2}}+\frac{m-2}{\sqrt{3}}+1-\frac{1}{\sqrt{3}}\\
&&+\left(\frac{m+1}{\sqrt{m+4}}-\frac{m-1}{\sqrt{m+3}}-\frac{1}{\sqrt{m+2}}\right)\\
&=&f(m)
\end{eqnarray*}
with equalities if and only if $d_G(v)=2$, $d_G(w)=5$, $G'=U_{2m,m}$
and $m=3$, i.e., $G=B_{6,3}$.

By combining Cases 1 and 2, the result follows.
\end{Proof}

\begin{Lemma}\label{l14}
Let $G\in\mathcal {B}(6,3)$. Then
$\chi(G)\ge\frac{4}{\sqrt{7}}+\frac{1}{\sqrt{6}}+1$ with equality if
and only if $G=B_{6,3}$.
\end{Lemma}
\begin{Proof}
If $G$ has a pendent vertex whose neighbor is of degree two, then
$G$ is the graph obtained from the unique $4$-vertex bicyclic graph
by attaching a path on two vertices to either a vertex of degree
three, or a vertex of degree two, and thus it is easily seen that
$\chi(G)>\frac{4}{\sqrt{7}}+\frac{1}{\sqrt{6}}+1$. Otherwise, by
Lemma \ref{l8},  $B_{6,3}$ is the unique graph with the minimum
sum--connectivity index.
\end{Proof}

Now we consider the bicyclic graphs with perfect matching. There is
a unique bicyclic graph with four vertices, and its matching number
is two.

\begin{Theorem}\label{th1}
Let $G\in\mathcal {B}(2m,m)$, where $m\ge3$. Then
\[
\chi(G)\ge\frac{m+1}{\sqrt{m+4}}+\frac{1}{\sqrt{m+3}}+\frac{m-3}{\sqrt{3}}+1
\]
with equality if and only if $G=B_{2m,m}$.
\end{Theorem}

\begin{Proof} Let
$f(m)=\frac{m+1}{\sqrt{m+4}}+\frac{1}{\sqrt{m+3}}+\frac{m-3}{\sqrt{3}}+1$.
We prove the result by induction on $m$. If $m=3$, then the result
follows from Lemma \ref{l14}.

Suppose that $m\ge4$ and the result holds for graphs in $\mathcal
{B}(2m-2,m-1)$. Let $G\in\mathcal {B}(2m,m)$ with a perfect matching
$M$.

If there is no pendent vertex with neighbor of degree two in $G$,
then by Lemma \ref{l8}, $\chi(G)>f(m)$. Suppose that $G$ has a
pendent vertex $u$ whose neighbor $v$ is of degree two. Then $uv\in
M$ and $G-u-v\in\mathcal {B}(2m-2,m-1)$. Let $w$ be the neighbor of
$v$ different from $u$. Since $|M|=m$, we have $d_G(w)\le m+2$. Note
that there is at most one pendent neighbor of $w$ in $G$. Then by Lemma
\ref{l5} (i), Lemma \ref{l6} (i) and the induction hypothesis,
\begin{eqnarray*}
\chi(G)&\ge&\chi(G-u-v)+
\frac{d_G(w)-1}{\sqrt{d_G(w)+2}}-\frac{d_G(w)-3}{\sqrt{d_G(w)+1}}-\frac{1}{\sqrt{d_G(w)}}+\frac{1}{\sqrt{3}}\\
&\ge&f(m-1)+
\frac{(m+2)-1}{\sqrt{(m+2)+2}}-\frac{(m+2)-3}{\sqrt{(m+2)+1}}-\frac{1}{\sqrt{m+2}}+\frac{1}{\sqrt{3}}\\
&=&f(m)
\end{eqnarray*}
with equalities if and only if $G-u-v=B_{2m-2,m-1}$ and
$d_G(w)=m+2$, i.e., $G=B_{2m,m}$.
\end{Proof}

In the following we consider the sum--connectivity indices of graphs
in the set of bicyclic graphs with $n$ vertices and matching number
$m$. We first consider the case $m\ge 3$.

\begin{Lemma} \cite{ZLW}\label{l4}
Let $G\in\mathcal {B}(n,m)$ with $n>2m \ge 6$, and $G$ has at least
one pendent vertex. Then there is a maximum matching $M$ and a
pendent vertex $u$ such that $u$ is not $M$-saturated.
\end{Lemma}

\begin{Theorem}\label{th2}
Let $G\in\mathcal {B}(n,m)$, where $3\le m \le\lfloor n/2 \rfloor$.
Then
\[
\chi(G)\ge
\frac{m+1}{\sqrt{n-m+4}}+\frac{n-2m+1}{\sqrt{n-m+3}}+\frac{m-3}{\sqrt{3}}+1
\]
with equality if and only if $G=B_{n,m}$.
\end{Theorem}

\begin{Proof}
Let
$f(n,m)=\frac{m+1}{\sqrt{n-m+4}}+\frac{n-2m+1}{\sqrt{n-m+3}}+\frac{m-3}{\sqrt{3}}+1$.
We prove the result by induction on $n$. If $n=2m$, then the result
follows from Theorem \ref{th1}. Suppose that $n>2m$ and the result
holds for graphs in $\mathcal {B}(n-1,m)$. Let $G\in\mathcal
{B}(n,m)$.

Suppose that there is no pendent vertex in $G$. Then $G\in
\widetilde{\mathbb{B}}(n)$ and $n=2m+1$. It is easily seen that
there are exactly three values for $\chi(G)$, and thus we have
$\chi(G)\ge \chi(H)=m-1+\frac{4}{\sqrt{6}}$ with $H\in {\bf
B}_2(2m+1)$. Let
$g(m)=\left(m-1+\frac{4}{\sqrt{6}}\right)-f(2m+1,m)=
\left(m-1+\frac{4}{\sqrt{6}}\right)
-\left(\frac{m+1}{\sqrt{m+5}}+\frac{2}{\sqrt{m+4}}+\frac{m-3}{\sqrt{3}}+1\right)$
for $m\ge 3$. Then
$g''(m)=(\frac{1}{4}m+\frac{17}{4})(m+5)^{-5/2}-\frac{3}{2}(m+4)^{-5/2}>0$,
and thus $g'(m)\ge g'(3)>0$, implying that $g(m)\ge g(3)>0$, i.e.,
$m-1+\frac{4}{\sqrt{6}}> f(2m+1,m)$. Then $\chi(G)>f(2m+1,m)$.

Suppose that there is at least one pendent vertex in $G$. By Lemma
\ref{l4}, there is a maximum matching $M$ and a pendent vertex $u$
of $G$ such that $u$ is not $M$-saturated. Then $G-u\in\mathcal
{B}(n-1,m)$. Let $v$ be the unique neighbor of $u$. Since $M$ is a
maximum matching, $M$ contains one edge incident with $v$. Note that
there are $n+1-m$ edges of $G$ outside $M$. Then $d_G(v)-1\le
n+1-m$, i.e., $d_G(v)\le n-m+2$. Let $s$ be the number of pendent
neighbors of $v$ in $G$.  Since at least $s-1$ pendent neighbors of
$v$ are not $M$-saturated, we have $s-1\le n-2m$, i.e.,  $s\le
n-2m+1$. By Lemma \ref{l5} (ii) with $k=n-2m+1$, Lemma \ref{l6} (ii)
and the induction hypothesis,
\begin{eqnarray*}
\chi(G)&\ge&\chi(G-u)+
\frac{d_G(v)-(n-2m+1)}{\sqrt{d_G(v)+2}}\\
&&+\frac{2(n-2m+1)-d_G(v)}{\sqrt{d_G(v)+1}}-\frac{(n-2m+1)-1}{\sqrt{d_G(v)}}\\
&\ge&f(n-1,m)+\frac{(n-m+2)-(n-2m+1)}{\sqrt{(n-m+2)+2}}\\
&&+\frac{2(n-2m+1)-(n-m+2)}{\sqrt{(n-m+2)+1}}-\frac{(n-2m+1)-1}{\sqrt{n-m+2}}\\
&=&f(n,m)
\end{eqnarray*}
with equalities if and only if $G-u=B_{n-1,m}$, $s=n-2m+1$ and
$d_G(v)=n-m+2$, i.e., $G=B_{n,m}$.
\end{Proof}

Now we  consider the sum--connectivity indices of graphs  bicyclic graphs  matching number
two. Let $B_n(a,b)$ be the graph obtained by attaching $a-3$ and
$b-3$ pendent vertices to the two vertices of degree three of the
unique $4$-vertex bicyclic  graph, respectively, where $a\ge b\ge
3$, $a+b=n+2$ and $n\ge4$.

\begin{Lemma}\label{l10}
Among the graphs in $\mathcal {B}(n,2)$ with $n\ge6$, $B_n(n-1,3)$
and $B_n(n-2,4)$ are respectively the unique graphs with the minimum
and the second minimum sum--connectivity indices, which are equal to
$\frac{1}{\sqrt{n+2}}+\frac{n-4}{\sqrt{n}}+\frac{2}{\sqrt{n+1}}+\frac{2}{\sqrt{5}}$
and
$\frac{1}{\sqrt{n+2}}+\frac{2}{\sqrt{n}}+\frac{n-5}{\sqrt{n-1}}+\frac{2}{\sqrt{6}}+\frac{1}{\sqrt{5}}$,
respectively.
\end{Lemma}
\begin{Proof}
Let $G\in\mathcal {B}(n,2)$. Then $G$ may be of three types:

(a) $G=B_n(a,b)$ with $a\ge b\ge 3$. Suppose that $a\ge b\ge 4$. Let
$f(x)=(x-4)x^{-1/2}+2(x+1)^{-1/2}$ for $x\ge 3$. Then
$f''(x)=-(\frac{1}{4}x+3)x^{-5/2}+\frac{3}{2}(x+1)^{-5/2}<0$,
implying that $f(x+1)-f(x)$ is decreasing for $x\ge 3$. It is easily
seen that
\begin{eqnarray*}
&&\chi(B_n(a+1,b-1))-\chi(B_n(a,b))\\
&=&[\chi(B_n(a+1,b-1))-\chi(B_{n-1}(a,b-1))]\\
&&-[\chi(B_n(a,b))-\chi(B_{n-1}(a,b-1))]\\
&=&\left(\frac{a-4}{\sqrt{a+2}}-\frac{a-3}{\sqrt{a+1}}+\frac{2}{\sqrt{a+3}}\right)
-\left(\frac{b-5}{\sqrt{b+1}}-\frac{b-4}{\sqrt{b}}+\frac{2}{\sqrt{b+2}}\right)\\
&=&[f(a+2)-f(a+1)]-[f(b+1)-f(b)]<0,
\end{eqnarray*}
and thus, $\chi(B_n(a,b))>\chi(B_n(a+1,b-1))$ for $a\ge b\ge 4$. It
follows that $B_n(n-1,3)$ and $B_n(n-2,4)$ are respectively the
unique graphs with the minimum and the second minimum
sum--connectivity indices, which are equal to
$\frac{1}{\sqrt{n+2}}+\frac{2}{\sqrt{n+1}}+\frac{n-4}{\sqrt{n}}+\frac{2}{\sqrt{5}}$
and
$\frac{1}{\sqrt{n+2}}+\frac{2}{\sqrt{n}}+\frac{n-5}{\sqrt{n-1}}+\frac{2}{\sqrt{6}}+\frac{1}{\sqrt{5}}$,
respectively.

(b) $G$ is the graph obtained by attaching $n-4$ pendent vertices to
a vertex of degree two of the unique $4$-vertex bicyclic graph. Then
\begin{eqnarray*}
&&\chi(G)=\frac{2}{\sqrt{n+1}}+\frac{n-4}{\sqrt{n-1}}+\frac{1}{\sqrt{6}}+\frac{2}{\sqrt{5}}\\
&>&\chi(B_n(n-2,4))=\frac{1}{\sqrt{n+2}}+\frac{2}{\sqrt{n}}
+\frac{n-5}{\sqrt{n-1}}+\frac{2}{\sqrt{6}}+\frac{1}{\sqrt{5}},
\end{eqnarray*}
since
$\chi(G)-\chi(B_n(n-2,4))=[g(n-1)-g(n)]+\frac{1}{\sqrt{5}}-\frac{1}{\sqrt{6}}>0$,
where
$g(x)=\frac{1}{\sqrt{x+2}}+\frac{1}{\sqrt{x}}-\frac{1}{\sqrt{x+1}}$
is decreasing for $x\ge 5$.

(c) $G$ is the graph obtained by attaching some pendent vertices to one or two
vertices of degree three of the unique $5$-vertex bicyclic graph
in  ${\bf B}^{(2)}_3(5)$, and by Lemma \ref{l7} and the arguments in
case (a), $\chi(G)>\chi(B_n(n-2,4))$.

Now the result follows easily.
\end{Proof}

\section{Minimum sum--connectivity index of bicyclic graphs}

In this section, we determine the minimum and the second minimum
sum--connectivity indices of bicyclic graphs with $n\ge 5$ vertices.

\begin{Theorem}
Among the  graphs in $\mathbb{B}(n)$ with $n\ge 5$, $B_n(n-1,3)$
is the unique graph with the minimum sum--connectivity index, which
is equal to
$\frac{1}{\sqrt{n+2}}+\frac{n-4}{\sqrt{n}}+\frac{2}{\sqrt{n+1}}+\frac{2}{\sqrt{5}}$,
the graph obtained by attaching a pendent vertex to a vertex of
degree two of the unique $4$-vertex bicyclic graph for $n=5$ is the
unique graph with the second minimum sum--connectivity index, which
is equal to $\frac{3}{\sqrt{6}}+\frac{2}{\sqrt{5}}+\frac{1}{2}$,
$B_n(n-2,4)$ for $n=6,7$ is the unique graph with the second minimum
sum--connectivity index, which is equal to
$\frac{1}{\sqrt{n+2}}+\frac{2}{\sqrt{n}}+\frac{n-5}{\sqrt{n-1}}+\frac{2}{\sqrt{6}}+\frac{1}{\sqrt{5}}$,
and  $B_{n,3}$ for $n\ge 8$ is the unique graph with the second
minimum sum--connectivity index, which is equal to
$\frac{4}{\sqrt{n+1}}+\frac{n-5}{\sqrt{n}}+1$.
\end{Theorem}

\begin{Proof} There are five graphs in $\mathbb{B}(5)$. Thus, the case $n=5$ may be checked
directly. Suppose in the following that $n\ge 6$.

Let $G\in \mathbb{B}(n)$ and $m$ the matching number of $G$, where
$2\le m \le\lfloor n/2 \rfloor$. If $m=2$, then by Lemma \ref{l10},
$\chi(G)\ge\chi(B_n(n-1,3))$ with equality if and only if
$G=B_n(n-1,3)$. If $m=3$, then by Theorem \ref{th2},
$\chi(G)\ge\chi(B_{n,3})$ with equality if and only if $G=B_{n,3}$.
If $m\ge4$, then by Theorem \ref{th2} and Lemma \ref{l7},
$\chi(G)\ge\chi(B_{n,m})>\chi(B_{n,m-1})>\cdots>\chi(B_{n,3})$. Let
$f(x)=\frac{1}{\sqrt{x}}-\frac{1}{\sqrt{x+1}}$ for $x\ge 6$. Then
$f''(x)=\frac{3}{4}x^{-5/2}-\frac{3}{4}(x+1)^{-5/2}>0$, implying
that $f(x+1)-f(x)$ is increasing for $x\ge 6$. Note that
\begin{eqnarray*}
&&\chi(B_{n,3})-\chi(B_n(n-1,3))\\
&=&\left(\frac{4}{\sqrt{n+1}}+\frac{n-5}{\sqrt{n}}+1\right)
-\left(\frac{1}{\sqrt{n+2}}+\frac{n-4}{\sqrt{n}}+\frac{2}{\sqrt{n+1}}+\frac{2}{\sqrt{5}}\right)\\
&=&f(n+1)-f(n)+1-\frac{2}{\sqrt{5}}\\
&\ge&f(7)-f(6)+1-\frac{2}{\sqrt{5}}>0.
\end{eqnarray*}
Thus $B_n(n-1,3)$ is the unique graph with the minimum
sum--connectivity index.

Suppose that $G\ne B_n(n-1,3)$. If $m=2$, then by Lemma \ref{l10}, $\chi(G)\ge \chi
(B_n(n-2,4))$ with equality if and only if $G=B_n(n-2,4)$. By the
arguments as above, the second minimum sum--connectivity index of
graphs in $\mathbb{B}(n)$ is precisely achieved by the minimum one
of $\chi(B_{n,3})$ and $\chi(B_n(n-2,4))$. If $n=6,7$, then
$\chi(B_{n,3})>\chi(B_n(n-2,4))$. Suppose that $n\ge 8$. Let
$g(x)=\frac{1}{\sqrt{x+1}}-\frac{3}{\sqrt{x}}-\frac{x-5}{\sqrt{x-1}}$
for $x\ge8$. Then
$g''(x)=\frac{3}{4}(x+1)^{-5/2}+\left[(\frac{1}{4}x+\frac{11}{4})(x-1)^{-5/2}-\frac{9}{4}x^{-5/2}\right]>0$,
implying that $g(x)-g(x+1)$ is decreasing for $x\ge 8$. Note that
\begin{eqnarray*}
&&\chi(B_{n,3})-\chi(B_n(n-2,4))\\
&=&\left(\frac{4}{\sqrt{n+1}}+\frac{n-5}{\sqrt{n}}+1\right)
-\left(\frac{1}{\sqrt{n+2}}+\frac{2}{\sqrt{n}}+\frac{n-5}{\sqrt{n-1}}+\frac{2}{\sqrt{6}}+\frac{1}{\sqrt{5}}\right)\\
&=&-\frac{1}{\sqrt{n+2}}+\frac{4}{\sqrt{n+1}}+\frac{n-7}{\sqrt{n}}
-\frac{n-5}{\sqrt{n-1}}+1-\frac{2}{\sqrt{6}}-\frac{1}{\sqrt{5}}\\
&=&g(n)-g(n+1)+1-\frac{2}{\sqrt{6}}-\frac{1}{\sqrt{5}}\\
&\le&g(8)-g(9)+1-\frac{2}{\sqrt{6}}-\frac{1}{\sqrt{5}}<0,
\end{eqnarray*}
and then $\chi(B_{n,3})<\chi(B_n(n-2,4))$. Thus
$B_n(n-2,4)$ for $n=6,7$
and  $B_{n,3}$ for $n\ge 8$ are the unique graphs with the second
minimum sum--connectivity index among graphs in $\mathbb{B}(n)$.
\end{Proof}

\section{Maximum sum--connectivity index of bicyclic graphs}

In this section, we determine the maximum and the second maximum
sum--connectivity indices of bicyclic graphs with $n\ge 5$ vertices.
Let $P_n$ be the path on $n$ vertices.

\begin{Lemma}\label{l1} \cite{ZhTr}
For a connected graph $Q$ with at least two vertices and a vertex
$u\in V(Q)$, let $G_1$ be the graph obtained from $Q$ by attaching
two paths $P_a$ and $P_b$ to $u$, $G_2$ the graph obtained from $Q$
by attaching a path $P_{a+b}$ to $u$, where $a\ge b\ge 1$. Then
$\chi(G_1)<\chi(G_2)$.
\end{Lemma}

\begin{Lemma}\label{l2}

Suppose that $M$ is a connected graph with $u\in V(M)$ and $2\le
d_M(u)\le 4$. Let $H$ be the graph obtained from $M$ by attaching a
path $P_a$ to $u$. Denote by $u_1$ and $u_2$ the two neighbors of
$u$ in $M$, and $u'$ the pendent vertex of the path attached to $u$
in $H$. Let $H'=H-uu_2+u'u_2$.
\begin{enumerate}
\item[(i)]
If $d_M(u)=2$ and the maximum degree of $M$ is at most five, then
$\chi(H')>\chi(H)$.

\item[(ii)]
If $d_M(u)=3$, and there are at least two neighbors of $u$ in $M$
with degree two and $d_M(u_2)=2$, then $\chi(H')>\chi(H)$.

\item[(iii)]
If $d_M(u)=4$ and all the neighbors of $u$ in $M$ are of degree two,
then $\chi(H')>\chi(H)$.

\end{enumerate}

\end{Lemma}

\begin{Proof}
(i) If $a=1$, then
\begin{eqnarray*}
&&\chi(H')-\chi(H)\\
&=&\left(\frac{1}{\sqrt{d_M(u_1)+2}}+\frac{1}{\sqrt{d_M(u_2)+2}}\right)
-\left(\frac{1}{\sqrt{d_M(u_1)+3}}+\frac{1}{\sqrt{d_M(u_2)+3}}\right)\\
&>&0.
\end{eqnarray*}
If $a\ge 2$, then
\begin{eqnarray*}
&&\chi(H')-\chi(H)\\
&=&\left(\frac{1}{\sqrt{d_M(u_1)+2}}-\frac{1}{\sqrt{d_M(u_1)+3}}\right)
+\left(\frac{1}{\sqrt{d_M(u_2)+2}}-\frac{1}{\sqrt{d_M(u_2)+3}}\right)\\
&&+1-\frac{1}{\sqrt{3}}-\frac{1}{\sqrt{5}}\\
&\ge&\left(\frac{1}{\sqrt{5+2}}-\frac{1}{\sqrt{5+3}}\right)
+\left(\frac{1}{\sqrt{5+2}}-\frac{1}{\sqrt{5+3}}\right)+1-\frac{1}{\sqrt{3}}-\frac{1}{\sqrt{5}}>0.
\end{eqnarray*}

(ii) There are two neighbors of $u$ with degree two,
let $d_1$ be the degree of the third neighbor of $u$ in $M$. If $a=1$, then since
$\frac{1}{2}+\frac{1}{\sqrt{5}}-\frac{2}{\sqrt{6}}>0$, we have
\begin{eqnarray*}
&&\chi(H')-\chi(H)\\
&=&\left(\frac{1}{\sqrt{d_1+3}}+\frac{1}{2}+\frac{2}{\sqrt{5}}\right)
-\left(\frac{1}{\sqrt{d_1+4}}+\frac{2}{\sqrt{6}}+\frac{1}{\sqrt{5}}\right)\\
&=&\left(\frac{1}{\sqrt{d_1+3}}-\frac{1}{\sqrt{d_1+4}}\right)
+\frac{1}{2}+\frac{1}{\sqrt{5}}-\frac{2}{\sqrt{6}}>0.
\end{eqnarray*}
If $a\ge 2$, then since
$1+\frac{2}{\sqrt{5}}-\frac{3}{\sqrt{6}}-\frac{1}{\sqrt{3}}>0$, we have
\begin{eqnarray*}
&&\chi(H')-\chi(H)\\
&=&\left(\frac{1}{\sqrt{d_1+3}}+1+\frac{2}{\sqrt{5}}\right)
-\left(\frac{1}{\sqrt{d_1+4}}+\frac{3}{\sqrt{6}}+\frac{1}{\sqrt{3}}\right)\\
&=&\left(\frac{1}{\sqrt{d_1+3}}-\frac{1}{\sqrt{d_1+4}}\right)+1
+\frac{2}{\sqrt{5}}-\frac{3}{\sqrt{6}}-\frac{1}{\sqrt{3}}>0.
\end{eqnarray*}

(iii) If $a=1$, then
\[
\chi(H')-\chi(H)=\left(\frac{1}{2}+\frac{4}{\sqrt{6}}\right)
-\left(\frac{4}{\sqrt{7}}+\frac{1}{\sqrt{6}}\right)>0.
\]
If $a\ge2$, then
\[
\chi(H')-\chi(H)=\left(1+\frac{4}{\sqrt{6}}\right)-\left(\frac{5}{\sqrt{7}}+\frac{1}{\sqrt{3}}\right)>0.
\]

The proof is now completed.
\end{Proof}

Let $\mathbb{B}_1(n)$ be the set of connected graphs on $n\ge 6$
vertices with exactly two vertex--disjoint cycles. Let
$\mathbb{B}_2(n)$ be the set of connected graphs on $n\ge 5$ vertices
with exactly two cycles of a common vertex. Let $\mathbb{B}_3(n)$
be the set of connected graphs on $n\ge 4$ vertices with exactly two
cycles with at least one edge in common. Obviously,
$\mathbb{B}(n)=\mathbb{B}_1(n)\cup\mathbb{B}_2(n)\cup\mathbb{B}_3(n)$.
For $u,v\in V(G)$, let $d_G(u,v)$ be the distance between $u$ and
$v$ in $G$.

\begin{Lemma}\label{l15}
Among the graphs in $\mathbb{B}_1(n)$ with $n\ge7$, the graphs in
${\bf B}_1^{(1)}(n)$ and the graphs in ${\bf B}_1^{(2)}(n)$ are
respectively the unique graphs with the maximum and the second
maximum sum--connectivity indices, which are equal to
$\frac{n-4}{2}+\frac{1}{\sqrt{6}}+\frac{4}{\sqrt{5}}$ and
$\frac{n-5}{2}+\frac{6}{\sqrt{5}}$, respectively.
\end{Lemma}

\begin{Proof}
Suppose that $G$ is a graph in $\mathbb{B}_1(n)\setminus\left\{{\bf
B}_1^{(1)}(n)\right\}$ with the maximum sum--connectivity index, and
$C^{(1)}$ and $C^{(2)}$ are its two cycles. Let $x_1\in V\left(C^{(1)}\right)$ and $y_1\in
V\left(C^{(2)}\right)$ be the two vertices such that
$d_G(x_1,y_1)=\min\{d_G(x,y):x\in V\left(C^{(1)}\right),y\in V\left(C^{(2)}\right)\}$. Let $Q$ be
the path joining $x_1$ and $y_1$. By Lemma \ref{l1},  the vertices
outside $C^{(1)}$, $C^{(2)}$ and $Q$ are of degree one or two, the vertices
on $C^{(1)}$, $C^{(2)}$ and $Q$ different from $x_1$ and $y_1$ are of degree
two or three, and $d_G(x_1),d_G(y_1)=3$ or $4$.

Suppose that $d_G(x_1,y_1)\ge 2$. If there is some vertex, say $x$,
on $C^{(1)}$, $C^{(2)}$ or $Q$ different from $x_1$ and $y_1$ with degree
three, then making use of Lemma \ref{l2} (i) to $H=G$ by setting
$u=x$, we may get a graph in $\mathbb{B}_1(n)\setminus\left\{{\bf
B}_1^{(1)}(n)\right\}$ with larger sum--connectivity index, a
contradiction. Thus the vertices on $C^{(1)}$, $C^{(2)}$ and $Q$ different
from $x_1$ and $y_1$ are of degree two. If $d_G(x_1)=4$, then making
use of Lemma \ref{l2} (ii) to $H=G$ by setting $u=x_1$, we may get a
graph in $\mathbb{B}_1(n)\setminus\left\{{\bf B}_1^{(1)}(n)\right\}$
with larger sum--connectivity index, a contradiction. Thus
$d_G(x_1)=3$. Similarly, we have $d_G(y_1)=3$. It follows that $G\in
{\bf B}_1^{(2)}(n)$.

Suppose that $d_G(x_1,y_1)=1$. Suppose that  one of $x_1$ and $y_1$,
say $x_1$, is of degree four.  Then by Lemma \ref{l2} (i), the
vertices on $C^{(1)}$ and $C^{(2)}$ different from $x_1$ and $y_1$ are of
degree two. If $d_G(y_1)=4$, then making use of Lemma \ref{l2} (ii)
to $H=G$ by setting $u=y_1$, we may get a graph in
$\mathbb{B}_1(n)\setminus\left\{{\bf B}_1^{(1)}(n)\right\}$ with
larger sum--connectivity index, a contradiction. Thus $d_G(y_1)=3$.
Denote by $x_2$ the pendent vertex of the path attached to $x_1$.
Consider $G_1=G-x_1y_1+x_2y_1\in {\bf B}_1^{(2)}(n)$. If
$d_G(x_1,x_2)=1$, then
\[
\chi(G_1)-\chi(G)=\frac{4}{\sqrt{5}}-\left(\frac{1}{\sqrt{7}}+\frac{2}{\sqrt{6}}+\frac{1}{\sqrt{5}}\right)>0.
\]
If $d_G(x_1,x_2)\ge 2$, then
\[
\chi(G_1)-\chi(G)=\left(\frac{1}{2}+\frac{4}{\sqrt{5}}\right)
-\left(\frac{1}{\sqrt{7}}+\frac{3}{\sqrt{6}}+\frac{1}{\sqrt{3}}\right)>0.
\]
In either case, $\chi(G_1)>\chi(G)$ with $G_1\in {\bf
B}_1^{(2)}(n)$, a contradiction. Thus $d_G(x_1)=d_G(y_1)=3$. Note
that $G\not\in{\bf B}_1^{(1)}(n)$ and by Lemma \ref{l2} (i), there
is exactly one vertex, say $x_3\in V\left(C^{(1)}\right)$, on $C^{(1)}$ and $C^{(2)}$
different from $x_1$ and $y_1$ with degree three. Denote by $x_4$
the pendent vertex of the path attached to $x_3$. Consider
$G_2=G-x_1y_1+x_4y_1\in {\bf B}_1^{(2)}(n)$. Let $d_1$ be the degree
of the neighbor of $x_4$, one neighbor of $x_1$ on $C^{(1)}$ is of degree two, and we denote by $d_2$ the other degree of the
neighbor of $x_1$ on $C^{(1)}$, where $d_1,d_2=2$ or $3$. We have
\begin{eqnarray*}
&&\chi(G_2)-\chi(G)\\
&=&\left(\frac{1}{\sqrt{d_1+2}}-\frac{1}{\sqrt{d_1+1}}\right)
+\left(\frac{1}{\sqrt{d_2+2}}-\frac{1}{\sqrt{d_2+3}}\right)+\frac{1}{2}-\frac{1}{\sqrt{6}}\\
&\ge&\left(\frac{1}{\sqrt{2+2}}-\frac{1}{\sqrt{2+1}}\right)
+\left(\frac{1}{\sqrt{3+2}}-\frac{1}{\sqrt{3+3}}\right)+\frac{1}{2}-\frac{1}{\sqrt{6}}>0,
\end{eqnarray*}
and thus, $\chi(G_2)>\chi(G)$ with $G_2\in {\bf B}_1^{(2)}(n)$,
which is also a contradiction.

Now we have shown that the graphs in ${\bf B}_1^{(2)}(n)$ are the
unique graphs in $\mathbb{B}_1(n)\setminus\left\{{\bf
B}_1^{(1)}(n)\right\}$ with the maximum sum--connectivity index.
Note that for $H_1\in{\bf B}_1^{(1)}(n)$ and $H_2\in {\bf
B}_1^{(2)}(n)$,
\[
\chi(H_1)=\frac{n-4}{2}+\frac{1}{\sqrt{6}}+\frac{4}{\sqrt{5}}>\chi(H_2)=\frac{n-5}{2}+\frac{6}{\sqrt{5}}.
\]
The result follows.
\end{Proof}

\begin{Lemma}\label{l16}
Among the graphs in $\mathbb{B}_3(n)$ with $n\ge5$, the graphs in
${\bf B}_3^{(1)}(n)$ and the graphs in ${\bf B}_3^{(2)}(n)$ are
respectively the unique graphs with the maximum and the second
maximum sum--connectivity indices, which are equal to
$\frac{n-4}{2}+\frac{1}{\sqrt{6}}+\frac{4}{\sqrt{5}}$ and
$\frac{n-5}{2}+\frac{6}{\sqrt{5}}$, respectively.
\end{Lemma}

\begin{Proof}
Suppose that $G$ is a graph in $\mathbb{B}_3(n)\setminus\left\{{\bf
B}_3^{(1)}(n)\right\}$ with the maximum sum--connectivity index.
Then $G$ has exactly three cycles, let $C^{(1)}$ and $C^{(2)}$ be its two
cycles such that the remaining one is of the maximum length. Let $A$
be the set of the common vertices of $C^{(1)}$ and $C^{(2)}$. Let $v_1$ and
$v_2$ be the two vertices in $A$ such that
$d_G(v_1,v_2)=\max\left\{d_G(x,y):x,y\in A\right\}$. By Lemma
\ref{l1}, the vertices outside $C^{(1)}$ and $C^{(2)}$ are of degree one or
two, the vertices on $C^{(1)}$ and $C^{(2)}$ different from $v_1$ and $v_2$
are of degree two or three, and $d_G(v_1),d_G(v_2)=3$ or $4$. Denote
by $v_1'$ ($v_2'$, respectively) the neighbor of $v_1$ on $C^{(1)}$
($v_2$ on $C^{(2)}$, respectively) different from the vertices in $A$.

If $d_G(v_1,v_2)\ge 2$, then by Lemma \ref{l2} (i) and (ii), we have $G\in {\bf B}_3^{(2)}(n)$.

Suppose that $d_G(v_1,v_2)=1$. Suppose that the lengths of $C^{(1)}$ and $C^{(2)}$ are at least four.
Consider $G_1=G-\{v_1v_1',v_2v_2'\}+\{v_1'v_2,v_1v_2'\}\in
\mathbb{B}_1(n)\setminus\left\{{\bf B}_1^{(1)}(n)\right\}$. Note
that
\begin{eqnarray*}
\chi(G_1)-\chi(G)
&=&\left(\frac{1}{\sqrt{d_G(v_1')+d_G(v_2)}}+\frac{1}{\sqrt{d_G(v_1)+d_G(v_2')}}\right)\\
&&-\left(\frac{1}{\sqrt{d_G(v_1)+d_G(v_1')}}+\frac{1}{\sqrt{d_G(v_2)+d_G(v_2')}}\right).
\end{eqnarray*}
If $d_G(v_1)=d_G(v_2)$, then $\chi(G_1)=\chi(G)$. If  $d_G(v_1)\ne
d_G(v_2)$, then by Lemma \ref{l2} (i), we have
$d_G(v_1')=d_G(v_2')=2$, and thus $\chi(G_1)=\chi(G)$. In either
case, we have $\chi(G_1)=\chi(G)$. By Lemma \ref{l15}, we have
$\chi(G)=\chi(G_1)\le \chi(H)=\frac{n-5}{2}+\frac{6}{\sqrt{5}}$ for
$H\in{\bf B}_1^{(2)}(n)$ with equality if and only if $G_1\in {\bf
B}_1^{(2)}(n)$, i.e., $G\in{\bf B}_3^{(2)}(n)$.

Suppose that at least one of $C^{(1)}$ and $C^{(2)}$, say $C^{(1)}$, is of
length three. Since $G\not\in {\bf B}_3^{(1)}(n)$, there are some vertices
outside $C^{(1)}$ and $C^{(2)}$. By Lemma \ref{l2} (i) and (ii), the
subgraph induced by the vertices outside $C^{(1)}$ and $C^{(2)}$ is a path,
say $P_k$, which is attached to $x\in V\left(C^{(1)}\right)\cup V\left(C^{(2)}\right)$. Suppose
that $x\ne v_1'$. Denote by $v_3$ the neighbor of $x$ outside $C^{(1)}$
and $C^{(2)}$. Consider
$G_2=G-xv_3+v_1'v_3\in\mathbb{B}_3(n)\setminus\left\{{\bf
B}_3^{(1)}(n)\right\}$.  If $x=v_1$ or $v_2$, then
\begin{eqnarray*}
&&\chi(G_2)-\chi(G)\\
&=&\left(\frac{1}{\sqrt{d_G(v_3)+3}}-\frac{1}{\sqrt{d_G(v_3)+4}}\right)
+\left(\frac{1}{\sqrt{6}}-\frac{1}{\sqrt{7}}\right)>0,
\end{eqnarray*}
and thus $\chi(G_2)>\chi(G)$, a contradiction. Hence $x\in
V\left(C^{(2)}\right)\setminus\{v_1,v_2\}$, and the length of $C^{(2)}$ is at least four.
Note that one  neighbor of $x$ on $C^{(2)}$ is of degree two. Denote
by $d_1$ the degree of the other neighbor of $x$ on $C^{(2)}$, where
$d_1=2$ or $3$. Then
\begin{eqnarray*}
&&\chi(G_2)-\chi(G)\\
&=&\left(\frac{1}{\sqrt{d_1+2}}-\frac{1}{\sqrt{d_1+3}}\right)+\frac{1}{2}+\frac{2}{\sqrt{6}}-\frac{3}{\sqrt{5}}\\
&\ge&\left(\frac{1}{\sqrt{3+2}}-\frac{1}{\sqrt{3+3}}\right)+\frac{1}{2}+\frac{2}{\sqrt{6}}-\frac{3}{\sqrt{5}}>0,
\end{eqnarray*}
and thus $\chi(G_2)>\chi(G)$, which is also a contradiction. Thus, $x=v_1'$. If
$k=1$, then
$\chi(G)=\frac{n-5}{2}+\frac{3}{\sqrt{6}}+\frac{2}{\sqrt{5}}+\frac{1}{2}$,
and if $k\ge 2$, then
$\chi(G)=\frac{n-6}{2}+\frac{3}{\sqrt{6}}+\frac{3}{\sqrt{5}}+\frac{1}{\sqrt{3}}$.
In either case, we have $\chi(G)< \frac{n-5}{2}+\frac{6}{\sqrt{5}}$.

Now we have shown that the graphs in ${\bf B}_3^{(2)}(n)$ are the
unique graphs in $\mathbb{B}_3(n)\setminus\left\{{\bf B}_3^{(1)}(n)\right\}$
with the maximum sum--connectivity index. Note that for $H_1\in{\bf
B}_3^{(1)}(n)$ and $H_2\in {\bf B}_3^{(2)}(n)$,
\[
\chi(H_1)=\frac{n-4}{2}+\frac{1}{\sqrt{6}}+\frac{4}{\sqrt{5}}>\chi(H_2)=\frac{n-5}{2}+\frac{6}{\sqrt{5}}.
\]
The result follows.
\end{Proof}

\begin{Theorem}
Among the graphs in $\mathbb{B}(n)$ with $n\ge 5$, the graph in
${\bf B}_3^{(1)}(5)$ and the graph in ${\bf B}_3^{(2)}(5)$ for $n=5$
are respectively the unique graphs with the maximum and the second
maximum sum--connectivity indices, the graphs in ${\bf
B}_1^{(1)}(6)\cup {\bf B}_3^{(1)}(6)$ and the graph in ${\bf
B}_3^{(2)}(6)$ for $n=6$ are respectively the unique graphs with the
maximum and the second maximum sum--connectivity indices, the graphs
in ${\bf B}_1^{(1)}(n)\cup {\bf B}_3^{(1)}(n)$ and the graphs in
${\bf B}_1^{(2)}(n)\cup {\bf B}_3^{(2)}(n)$ for $n\ge 7$ are
respectively the unique graphs with the maximum and the second
maximum sum--connectivity indices, where
$\chi(G)=\frac{n-4}{2}+\frac{1}{\sqrt{6}}+\frac{4}{\sqrt{5}}$ for
$G\in {\bf B}_1^{(1)}(n)\cup {\bf B}_3^{(1)}(n)$ and
$\chi(H)=\frac{n-5}{2}+\frac{6}{\sqrt{5}}$ for $H\in {\bf
B}_1^{(2)}(n)\cup {\bf B}_3^{(2)}(n)$.
\end{Theorem}

\begin{Proof}
Suppose that $G$ is a graph in $\mathbb{B}_2(n)$ with the maximum
sum--connectivity index, and $C^{(1)}$ and $C^{(2)}$ are its two cycles. Let
$u$ be the unique common vertex of $C^{(1)}$ and $C^{(2)}$. By Lemma
\ref{l1}, the vertices outside $C^{(1)}$ and $C^{(2)}$ are of degree one or
two, the vertices on $C^{(1)}$ and $C^{(2)}$ different from $u$ are of
degree two or three, and $d_G(u)=4$ or $5$. Moreover, by Lemma
\ref{l2} (i), the vertices on $C^{(1)}$ and $C^{(2)}$ different from $u$ are
of degree two. If $d_G(u)=5$, then making use of Lemma \ref{l2}
(iii) to $H=G$, we may get a graph in
$\mathbb{B}_2(n)$ with larger sum--connectivity index, a
contradiction. Thus $d_G(u)=4$, i.e., $G\in {\bf B}_2(n)$.

Note that for $H_1\in{\bf B}_1^{(1)}(n)$, $H_1'\in{\bf
B}_1^{(2)}(n)$, $H_2\in {\bf B}_2(n)$, $H_3\in{\bf B}_3^{(1)}(n)$
and $H_3'\in{\bf B}_3^{(2)}(n)$,
\begin{eqnarray*}
&&\chi(H_1)=\chi(H_3)=\frac{n-4}{2}+\frac{1}{\sqrt{6}}+\frac{4}{\sqrt{5}}\\
&>&\chi(H_1')=\chi(H_3')=\frac{n-5}{2}+\frac{6}{\sqrt{5}}\\
&>&\chi(H_2)=\frac{n-3}{2}+\frac{4}{\sqrt{6}}.
\end{eqnarray*}
Then the result follows from Lemmas \ref{l15} and \ref{l16}.
\end{Proof}

\vspace{6mm}

\noindent {\it Acknowledgement.\/} This work was supported by the
National Natural Science Foundation of China (Grant No.~10671076).

\end{document}